\documentclass[12pt,a4]{amsart}

\usepackage[mathscr]{eucal}
\usepackage{amsmath}
\usepackage{graphicx}
\usepackage{amssymb}
\usepackage{latexsym}
\usepackage{amsthm}
\usepackage{amscd}
\usepackage{enumerate}

\theoremstyle{plain}
\newtheorem{thm}{Theorem}[section]

\newtheorem{cor}[thm]{Corollary}

\theoremstyle{definition}
\newtheorem{Def}[thm]{Definition}

\numberwithin{equation}{section}

\title{A combinatorial proof of an identity for the divisor generating function
}
\author{Masanori Ando\\ (Okayama university)}
\date{}
\begin{document}
\pagestyle{empty}
\maketitle\thispagestyle{empty}

\section{Introduction}

In \cite{U}, 
Uchimura proved the following $q$-series identity
\begin{equation}
\label{eq;U}
\displaystyle \sum_{k \geq 1}{(-1)^k \frac {q^{\frac {k(k+1) }{2 }
}
 }{(q;q )_k (1-q^k) } }
=\sum_{k \geq 1 }{\frac {q^k }{1-q^k } },  
\end{equation}
where $(a;b)_k=(1-a)(1-ab)(1-ab^2)\cdots (1-ab^{k-1})$. 
The formula \eqref{eq;U} has been known since 100 years ago (\cite{K}). 
This identity is an infinite version of the following $q$-series identity called {\it {Problem 6407}} in American Mathematical Monthly\cite{vh}. 
\[
\sum_{k=1}^{m}{(-1)^{k-1}\left [ 
\begin{array}{l} 
m \\
k 
\end{array}
\right ] 
\frac{q^{\frac {k(k+1)}{2} }}{1-q^k} 
}
=\sum _{k=1}^{m} {\frac {q^k}{1-q^k}}. 
\]
Here $\left [ 
\begin{array}{l} 
m \\
k 
\end{array}
\right ]$ is a $q-$binomial coefficient. 
Many authors have generalized these identities (see e.g. \cite{vgcz}). 
In this paper, we translate these identities and Dilcher's generalization \cite{D} into combinatorics of partitions, and give a combinatorial proof of them. 
For example, we transform \eqref{eq;U} as 
\[
\sum_{\lambda \in {\mathcal {SP} } }{(-1)^{\ell (\lambda)-1 }\lambda _{\ell (\lambda) } q^{|\lambda | } } =\sum_{n \geq 1}{\sigma _0 (n)q^n }. 
\]
It is a $q$-series identity about {\it {strict partitions}} and a {\it {divisor function}}. 

\bigskip
The generalizations of \eqref{eq;U} we give in this paper are the following. 
\[
\displaystyle \sum_{k= 1}^{\infty }{(-1)^k \frac {b^kq^{\frac {k(k+1) }{2 }
+(m-1)k }
 }{(bq;q )_k (1-q^k)^m } }
=\sum_{j_1= 1 }^{\infty }{\frac {b^{j_1 }q^{j_1} }{1-q^{j_1} } 
\sum_{j_2 =1 }^{j_1 }{\frac {q^{j_2} }{1-q^{j_2} } 
\cdots
\sum_{j_m =1 }^{j_{m-1 } }{\frac {q^{j_m} }{1-q^{j_m} } 
}
}
}, 
\]
\[
\displaystyle \sum_{k= 1}^{t }{(-1)^k \frac {b^kq^{\frac {k(k+1) }{2 }
+(m-1)k }
 }{(1-q^k)^m } 
\left [ 
\begin{array}{l} 
t \\
k 
\end{array}
\right ] _{q,b}
}
=\sum_{j_1= 1 }^{t }{\frac {b^{j_1 }q^{j_1} }{1-q^{j_1} } 
\sum_{j_2 =1 }^{j_1 }{\frac {q^{j_2} }{1-q^{j_2} } 
\cdots
\sum_{j_m =1 }^{j_{m-1 } }{\frac {q^{j_m} }{1-q^{j_m} } 
}
}
}. 
\]
As a by-product of their proofs, we obtain some combinatorial results. 

\section{Young diagrams}

\begin{Def}
Let $n$ be a positive integer. A partition $\lambda $ of $n$ is an integer sequence 
\[
\lambda =(\lambda_1,\lambda_2,\ldots,\lambda_\ell) 
\]
satisfying $\lambda _1 \geq \lambda _2 \geq \ldots \geq \lambda _\ell >0$ and $\sum _{i=1}^{\ell }{\lambda _i} =n$. 
We call $\ell (\lambda):= \ell$ the length of $\lambda $, 
and each $\lambda _i$ a part of $\lambda $. 
The set of partitions of $n$ is denoted by ${\mathcal {P}}(n) $. 
\end{Def}
\begin{Def}
A partition $\lambda $ is said to be strict if $\lambda _1 > \lambda _2 > \ldots > \lambda _\ell >0$. 
The set of strict partitions of $n$ is denoted by ${\mathcal {SP}}(n) $. 
\end{Def}
\begin{Def}
Let $\lambda = (\lambda_1,\lambda_2,\ldots,\lambda_\ell) $ be a partition. The Young diagram of $\lambda $ is defined by
\[
Y(\lambda):=\{ (i,j)\in \mathbb{N}\times \mathbb{N} \ |\  1\leq i\leq \ell ,1\leq j \leq \lambda_i \}. 
\]
We call $(i,j)\in Y(\lambda )$ the (i,j)-cell of $\lambda$. 
And the set of the corners of $\lambda $ is defined by 
\[
C(\lambda ):= 
\{ (i,j)\in Y(\lambda )\ |\ (i+1,j ), (i, j+1 )\not \in Y(\lambda ) 
\}. 
\]
We put $c(\lambda ):=\sharp C(\lambda )$, the number of the corners of $\lambda $. 
\end{Def}
\begin{Def}
Let $(i,j)\in Y(\lambda)$, The $(i,j)$-hook length of $\lambda $ is defined by
\[
h_{ij}(\lambda):=\sharp \{ (i',j')\in Y(\lambda)\  |\  i'=i,j'\geq j\ \  \textrm{or}\ \  j'=j,i'\geq i \ \}
\]
And we put $a_{ij}(\lambda):= \lambda_i - j+1$, the $(i,j)$-arm length of $\lambda $. 
We remark that our definition of arm length $a_{ij}$ is different by 1 from the usual definition \cite{Mac}. 
\end{Def}

\section{$q$-series identity for the divisor function}

\begin{thm}$($Uchimura's identity$)$
\[
\displaystyle \sum_{k \geq 1}{(-1)^k \frac {q^{\frac {k(k+1) }{2 }
}
 }{(q;q )_k (1-q^k) } }
=\sum_{k \geq 1 }{\frac {q^k }{1-q^k } },  
\]
where $(a;b)_k=(1-a)(1-ab)(1-ab^2)\cdots (1-ab^{k-1})$. 
\end{thm}
Remark that the right-hand side is computed as 
\[
\sum _{k \geq 1}{\frac {q^k}{1-q^k}} = \sum _{k=1}^{\infty }{(q^k+q^{2k}+q^{3k}+ \cdots )}=\sum _{n \geq 1}{\sigma _0 (n)q^n}, 
\]
where $\sigma _0 (n)$ is the number of positive divisors of $n$. 
We now translate this identity into a language of Young diagrams. Then we are able to prove this identity combinatorially. \\
\newpage
\noindent
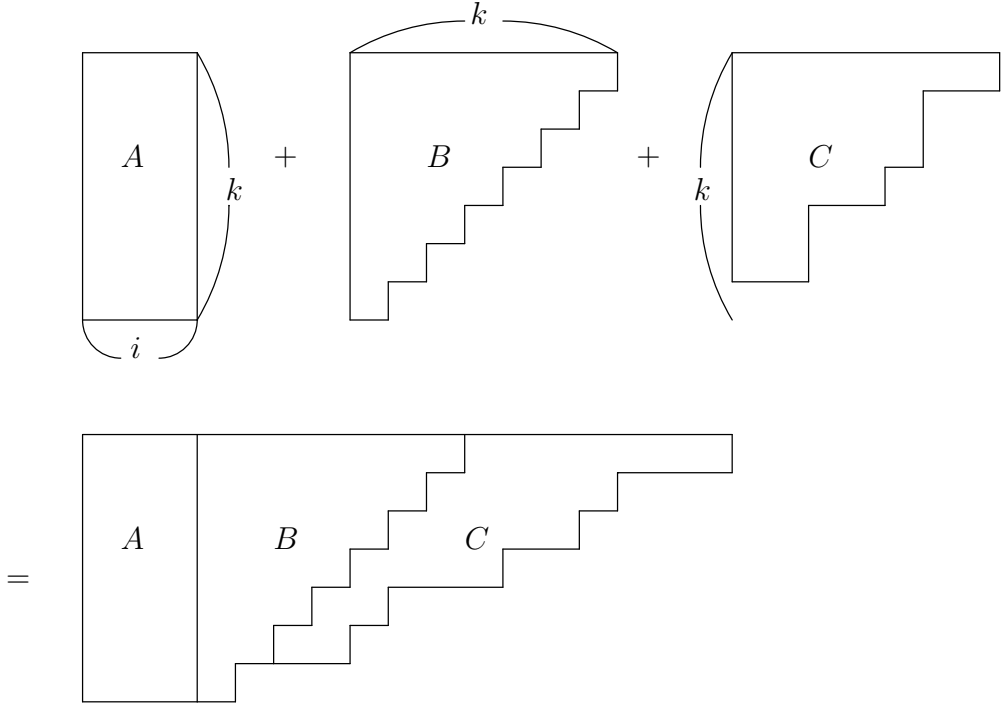
\textbf{Figure 1.}\\
\unitlength 0.1in
\begin{picture}(55.00,40.70)(3.00,-42.00)
%
\special{pn 8}%
\special{pa 600 600}%
\special{pa 1200 600}%
\special{fp}%
\special{pa 1200 600}%
\special{pa 1200 2000}%
\special{fp}%
\special{pa 1200 2000}%
\special{pa 600 2000}%
\special{fp}%
\special{pa 600 2000}%
\special{pa 600 600}%
\special{fp}%
%
\special{pn 8}%
\special{ar 200 1400 1166 1166  6.2831853 6.2831853}%
\special{ar 200 1400 1166 1166  0.0000000 0.5404195}%
%
\special{pn 8}%
\special{ar 200 1200 1166 1166  5.7427658 6.2831853}%
%
\special{pn 8}%
\special{ar 800 2000 200 200  1.5707963 3.1415927}%
%
\special{pn 8}%
\special{ar 1000 2000 200 200  6.2831853 6.2831853}%
\special{ar 1000 2000 200 200  0.0000000 1.5707963}%
%
\special{pn 8}%
\special{pa 2000 600}%
\special{pa 2000 2000}%
\special{fp}%
\special{pa 2000 2000}%
\special{pa 2200 2000}%
\special{fp}%
\special{pa 2200 2000}%
\special{pa 2200 1800}%
\special{fp}%
\special{pa 2200 1800}%
\special{pa 2400 1800}%
\special{fp}%
\special{pa 2400 1800}%
\special{pa 2400 1600}%
\special{fp}%
\special{pa 2400 1600}%
\special{pa 2600 1600}%
\special{fp}%
\special{pa 2600 1600}%
\special{pa 2600 1400}%
\special{fp}%
\special{pa 2600 1400}%
\special{pa 2800 1400}%
\special{fp}%
\special{pa 2800 1400}%
\special{pa 2800 1200}%
\special{fp}%
\special{pa 2800 1200}%
\special{pa 3000 1200}%
\special{fp}%
\special{pa 3000 1200}%
\special{pa 3000 1000}%
\special{fp}%
\special{pa 3000 1000}%
\special{pa 3200 1000}%
\special{fp}%
\special{pa 3200 1000}%
\special{pa 3200 800}%
\special{fp}%
\special{pa 3200 800}%
\special{pa 3400 800}%
\special{fp}%
\special{pa 3400 800}%
\special{pa 3400 600}%
\special{fp}%
\special{pa 3400 600}%
\special{pa 2000 600}%
\special{fp}%
%
\special{pn 8}%
\special{pa 4000 600}%
\special{pa 4000 1800}%
\special{fp}%
\special{pa 4000 1800}%
\special{pa 4400 1800}%
\special{fp}%
\special{pa 4400 1800}%
\special{pa 4400 1400}%
\special{fp}%
\special{pa 4400 1400}%
\special{pa 4800 1400}%
\special{fp}%
\special{pa 4800 1400}%
\special{pa 4800 1200}%
\special{fp}%
\special{pa 4800 1200}%
\special{pa 5000 1200}%
\special{fp}%
\special{pa 5000 1200}%
\special{pa 5000 800}%
\special{fp}%
\special{pa 5000 800}%
\special{pa 5400 800}%
\special{fp}%
\special{pa 5400 800}%
\special{pa 5400 600}%
\special{fp}%
\special{pa 5400 600}%
\special{pa 4000 600}%
\special{fp}%
%
\special{pn 8}%
\special{pa 600 2600}%
\special{pa 4000 2600}%
\special{fp}%
\special{pa 4000 2600}%
\special{pa 4000 2800}%
\special{fp}%
\special{pa 4000 2800}%
\special{pa 3400 2800}%
\special{fp}%
\special{pa 600 2600}%
\special{pa 600 4000}%
\special{fp}%
\special{pa 600 4000}%
\special{pa 1200 4000}%
\special{fp}%
\special{pa 1200 4000}%
\special{pa 1200 2600}%
\special{fp}%
\special{pa 1200 4000}%
\special{pa 1400 4000}%
\special{fp}%
\special{pa 1400 4000}%
\special{pa 1400 3800}%
\special{fp}%
\special{pa 1400 3800}%
\special{pa 1600 3800}%
\special{fp}%
\special{pa 1600 3800}%
\special{pa 1600 3600}%
\special{fp}%
\special{pa 1600 3600}%
\special{pa 1800 3600}%
\special{fp}%
\special{pa 1800 3600}%
\special{pa 1800 3400}%
\special{fp}%
\special{pa 1800 3400}%
\special{pa 2000 3400}%
\special{fp}%
\special{pa 2000 3400}%
\special{pa 2000 3200}%
\special{fp}%
\special{pa 2000 3200}%
\special{pa 2200 3200}%
\special{fp}%
\special{pa 2200 3200}%
\special{pa 2200 3000}%
\special{fp}%
\special{pa 2200 3000}%
\special{pa 2400 3000}%
\special{fp}%
\special{pa 2400 3000}%
\special{pa 2400 2800}%
\special{fp}%
\special{pa 2400 2800}%
\special{pa 2600 2800}%
\special{fp}%
\special{pa 2600 2800}%
\special{pa 2600 2600}%
\special{fp}%
\special{pa 1600 3800}%
\special{pa 2000 3800}%
\special{fp}%
\special{pa 2000 3800}%
\special{pa 2000 3600}%
\special{fp}%
\special{pa 2000 3600}%
\special{pa 2200 3600}%
\special{fp}%
\special{pa 2200 3600}%
\special{pa 2200 3400}%
\special{fp}%
\special{pa 2200 3400}%
\special{pa 2800 3400}%
\special{fp}%
\special{pa 2800 3400}%
\special{pa 2800 3200}%
\special{fp}%
\special{pa 2800 3200}%
\special{pa 3200 3200}%
\special{fp}%
\special{pa 3200 3200}%
\special{pa 3200 3000}%
\special{fp}%
\special{pa 3200 3000}%
\special{pa 3400 3000}%
\special{fp}%
\special{pa 3400 2800}%
\special{pa 3400 2800}%
\special{fp}%
%
\special{pn 8}%
\special{pa 3400 2800}%
\special{pa 3400 3000}%
\special{fp}%

\special{pn 8}%
\special{ar 5000 1200 1166 1166  3.1415927 3.6820122}%

\special{pn 8}%
\special{ar 5000 1400 1166 1166  2.6011732 3.1415927}%

\special{pn 8}%
\special{ar 2600 1600 1166 1166  4.1719695 4.7123890}%
%
\special{pn 8}%
\special{ar 2800 1600 1166 1166  4.7123890 5.2528085}%

\put(26.3000,-4.5000){\makebox(0,0)[lb]{$k$}}%

\put(13.5000,-13.7000){\makebox(0,0)[lb]{$k$}}%

\put(8.5000,-22.0000){\makebox(0,0)[lb]{$i$}}%

\put(16.0000,-12.0000){\makebox(0,0)[lb]{$+$}}%

\put(35.0000,-12.0000){\makebox(0,0)[lb]{$+$}}%

\put(38.0000,-13.7000){\makebox(0,0)[lb]{$k$}}%

\put(8.0000,-12.0000){\makebox(0,0)[lb]{$A$}}%

\put(24.0000,-12.0000){\makebox(0,0)[lb]{$B$}}%

\put(44.0000,-12.0000){\makebox(0,0)[lb]{$C$}}%

\put(8.0000,-32.0000){\makebox(0,0)[lb]{$A$}}%

\put(16.0000,-32.0000){\makebox(0,0)[lb]{$B$}}%

\put(26.0000,-32.0000){\makebox(0,0)[lb]{$C$}}%

\put(2.0000,-34.0000){\makebox(0,0)[lb]{$=$}}%
\end{picture}%
\\
Looking at each term of the left-hand side, 
$\displaystyle q^{\frac {k(k+1)}{2}}$ is translated into the stairs $B$ in Figure $1$. 
Since $\displaystyle \frac {1}{(q;q)_k}$ is the generating function of partitions whose lengths are at most $k$, 
this term corresponds to $C$ in Figure $1$. 
The leftover $\displaystyle \frac{1}{1-q^k}$ is the generating function of rectangular Young diagrams whose vertical lengths are equal to $k$. 
This part corresponds to $A$ in Figure $1$. 
Therefore the left-hand side of the identity is an alternating sum over $k \geq 1$ of $A+B+C$. 
As is noted in Figure $1$, the ``sum" $A+B+C$ is a strict partition.
For a strict partition $\lambda $, we count the number of tuples $(A,B,C)$ such that $A+B+C = \lambda$. 
Let $\lambda $ be a fixed strict partition of length $k$. 
One can embed the stairs $B$ into $\lambda $ in $\lambda _k$ ways. 
For each embedding the rectangle $A$ and the partition $C$ are uniquely determined, respectively. 
Therefore there are $\lambda _k$ tuples $(A,B,C)$, such that $A+B+C = \lambda $. 
Summing up over $k$, Theorem 3.1 reads 
\begin{eqnarray}
\label{id;div}
\sum_{\lambda \in {\mathcal {SP} } }{(-1)^{\ell (\lambda)-1 }\lambda _{\ell (\lambda) } q^{|\lambda | } } =\sum_{n \geq 1}{\sigma _0 (n)q^n }. 
\end{eqnarray}
\begin{thm}
For any positive integers $n$ and $k$, 
\begin{eqnarray*}
&&\sharp \{ \lambda \in {\mathcal {SP} }(n) \ |\  \lambda_1 \geq k > \lambda _1 -\lambda _{\ell (\lambda)},  \ell (\lambda ): \textrm {odd} \} \\
&-&\sharp \{ \lambda \in {\mathcal {SP} }(n) \ |\  \lambda_1 \geq k > \lambda _1 -\lambda _{\ell (\lambda)}, \ell (\lambda ): \textrm {even} \}\\
&=&\left \{ \begin{array}{ll}
1  &(k\mid n) \\
0  &(k\nmid n) .
\end{array} \right.
\end{eqnarray*}
\end{thm}
\noindent
\textbf{Example}. @For $n=5$, 
we draw Young diagrams $Y(\lambda)$ of all strict partitions of $5$, and write arm length in $(1,j)$-cell for $1\leq j \leq \lambda _{\ell (\lambda )}$. \\
\\
\textbf{Figure 2.}\\
\unitlength 0.1in
\begin{picture}(38.00,10.00)(4.00,-14.00)
%
\special{pn 8}%
\special{pa 400 400}%
\special{pa 400 400}%
\special{fp}%
%
\special{pn 8}%
\special{pa 600 600}%
\special{pa 1600 600}%
\special{fp}%
\special{pa 1600 600}%
\special{pa 1600 800}%
\special{fp}%
\special{pa 1600 800}%
\special{pa 600 800}%
\special{fp}%
\special{pa 600 800}%
\special{pa 600 600}%
\special{fp}%
%
\special{pn 8}%
\special{pa 2200 600}%
\special{pa 3000 600}%
\special{fp}%
\special{pa 3000 600}%
\special{pa 3000 800}%
\special{fp}%
\special{pa 3000 800}%
\special{pa 2200 800}%
\special{fp}%
\special{pa 2200 600}%
\special{pa 2200 1000}%
\special{fp}%
\special{pa 2200 1000}%
\special{pa 2400 1000}%
\special{fp}%
\special{pa 2400 1000}%
\special{pa 2400 600}%
\special{fp}%
\special{pa 2600 600}%
\special{pa 2600 800}%
\special{fp}%
\special{pa 2800 800}%
\special{pa 2800 600}%
\special{fp}%
\special{pa 1400 600}%
\special{pa 1400 800}%
\special{fp}%
\special{pa 1200 800}%
\special{pa 1200 600}%
\special{fp}%
\special{pa 1000 600}%
\special{pa 1000 800}%
\special{fp}%
\special{pa 800 800}%
\special{pa 800 600}%
\special{fp}%
%
\special{pn 8}%
\special{pa 3600 600}%
\special{pa 4200 600}%
\special{fp}%
\special{pa 4200 600}%
\special{pa 4200 800}%
\special{fp}%
\special{pa 4200 800}%
\special{pa 3600 800}%
\special{fp}%
\special{pa 3600 600}%
\special{pa 3600 1000}%
\special{fp}%
\special{pa 3600 1000}%
\special{pa 4000 1000}%
\special{fp}%
\special{pa 4000 1000}%
\special{pa 4000 600}%
\special{fp}%
\special{pa 3800 600}%
\special{pa 3800 1000}%
\special{fp}%

\put(6.4000,-7.5000){\makebox(0,0)[lb]{$5$}}%

\put(8.4000,-7.5000){\makebox(0,0)[lb]{$4$}}%

\put(10.4000,-7.5000){\makebox(0,0)[lb]{$3$}}%

\put(12.4000,-7.5000){\makebox(0,0)[lb]{$2$}}%

\put(14.4000,-7.5000){\makebox(0,0)[lb]{$1$}}%

\put(22.4000,-7.5000){\makebox(0,0)[lb]{$4$}}%

\put(18.2000,-7.7000){\makebox(0,0)[lb]{$-$}}%

\put(32.2000,-7.7000){\makebox(0,0)[lb]{$-$}}%

\put(36.4000,-7.5000){\makebox(0,0)[lb]{$3$}}%

\put(38.4000,-7.5000){\makebox(0,0)[lb]{$2$}}%

\put(4.0000,-14.0000){\makebox(0,0)[lb]{$=$}}%

\put(8.0000,-14.5000){\makebox(0,0)[lb]{$5,1$. }}%
\end{picture}%
\\
\\
Here the numbers are regarded as variables. 
Number $5$ and $1$ are the positive divisors of $5$. 

\bigskip
\noindent
\textbf{Proof of Theorem 3.2}\\
We consider the set of strict partitions of $n$ such that $a_{1,j}=k $ for $1\leq j \leq \lambda _{\ell (\lambda )}$:
\[
{\mathfrak{D}}(n,k):= \{ \lambda \in {\mathcal {SP}}(n) \ |\  \lambda _1 - \lambda _{\ell (\lambda )} <  k \leq \lambda _1 \}. 
\]
We divide these strict partitions into two classes $A$ and $B$: 
\[
A=\{ \lambda \in {\mathfrak{D}}(n,k)\ |\ {}^\forall i, \ k\nmid \lambda_i\}, \ \ \ B=\{ \lambda \in {\mathfrak{D}}(n,k)\ |\ {}^\exists  i, \ k\mid \lambda_i\}. 
\]
We consider a map between them that changes the length by $1$. 
\[
\alpha_k :\ A \rightarrow B , \ \ \ \ \ \alpha _k(\lambda )= \lambda ' , 
\]
where $\lambda ' \in B$ is defined in the following steps: 
\begin{enumerate}[Step 1.]
\item Append $0$ in the tail of $\lambda $ to get $(\lambda _1 , \ldots , \lambda _{\ell +1})$. 
\item Subtract $k$ from $\textrm {max} (\lambda _1 , \ldots , \lambda _{\ell+1})$, and add $k$ to $\lambda _{\ell+1}$. 
\item Repeat Step 2 till $\textrm {max} (\lambda _1 , \ldots , \lambda _{\ell+1}) - \textrm {min} (\lambda _1 , \ldots , \lambda _{\ell+1})$ 
gets less than $k$. 
\item From the resulting composition we have the partition $\lambda '= (\lambda _1, \ldots , \lambda _{\ell +1})$, by arranging parts. 
\end{enumerate}
By the above construction, we have $\ell (\lambda ')=\ell (\lambda )+1$, and
\[
\sharp \{i \ |\ \lambda _i  ({\textrm {mod}} k)\equiv j , 1\leq i \leq \ell \}
=\sharp \{i \ |\ \lambda '_i  ({\textrm {mod}} k)\equiv  j , 1\leq i \leq \ell \}
\]
for $1\leq j \leq k-1$. 
The partition that $\alpha _k$ can not pair up is $(\lambda_1 , \ldots , \lambda _\ell )\equiv (0) \ \ \textrm {mod} k$. 
Therefore $\lambda = (n) $ is left when $n$ is a multiple of $k$. 
\qed

\bigskip 
\noindent
\textbf{Proof of theorem 3.1}\\
Sum of the left-hand side of Theorem 3.2 over $k$ is 
\begin{eqnarray*}
&&\sum_{k \geq 1} \sharp \{ \lambda \in {\mathcal {SP} }(n) \ |\  \lambda_1 \geq k > \lambda _1 -\lambda _{\ell (\lambda )}, \ell (\lambda ): \textrm {odd} \} \\
&-&\sum _{k \geq 1}\sharp \{ \lambda \in {\mathcal {SP} }(n) \ |\  \lambda_1 \geq k > \lambda _1 -\lambda _{\ell (\lambda )}, \ell (\lambda ): \textrm {even} \}\\
=&&\sum _{\lambda \in {\mathcal {SP} }(n) }{(-1)^{\ell (\lambda -1)}\lambda _{\ell (\lambda )} }. 
\end{eqnarray*}
And sum of the right-hand side is 
\begin{eqnarray*}
\sum _{k\mid n}{1}
=\sigma _0 (n). 
\end{eqnarray*}
They are the coefficients of $q^n$ in \eqref{id;div}. 
\qed

\bigskip
Theorem 3.1 is the generating function for the total sum of Theorem 3.2. 
Taking the sum over $k$ from $1$ to $m$ for Theorem 3.1, 
we have the identity of Problem 6407 which is shown in the following generating function. 
\[
\sum_{k=1}^{m}{(-1)^{k-1}\left [ 
\begin{array}{l} 
m \\
k 
\end{array}
\right ] 
\frac{q^{\frac {k(k+1)}{2} }}{1-q^k} 
}
=\sum _{k=1}^{m} {\frac {q^k}{1-q^k}}, 
\]
where the $q$-binomial coefficient is defined by
\[
\left [ 
\begin{array}{l} 
m \\
k 
\end{array}
\right ]
=
\left \{
\begin{array}{cc}
\frac{(q;q)_n}{(q;q)_k(q;q)_{n-k}} & 1\leq k \leq n \\
0 & \textrm{otherwise} . 
\end{array}
\right.
\] 
\begin{cor}
For $n=2(2m+1)$, 
\begin{eqnarray*}
&&\sum_{\lambda \in {\mathcal {SP} }(n)}\sharp \{ h_{1,i}(\lambda )\ | \ i\leq \lambda_{\ell (\lambda)}, \ \textrm {odd} \} \\
&=&\sum_{\lambda \in {\mathcal {SP} }(n)}\sharp \{h_{1,i}(\lambda )\ | \ i\leq \lambda_{\ell (\lambda)}, \ \textrm {even} \}. 
\end{eqnarray*}
\end{cor}
\noindent
\textbf{Example.} For $n=6$, 
we draw Young diagrams $Y(\lambda)$ of all strict partitions of $6$, and write hook lengths in $(1,j)$-cell for $1\leq j \leq \lambda _{\ell (\lambda)}$. \\
\\
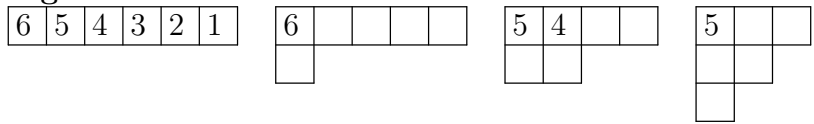
\textbf{Figure 3.}
\\
\unitlength 0.1in
\begin{picture}(42.00,6.00)(4.00,-10.00)
%
\special{pn 8}%
\special{pa 400 400}%
\special{pa 1600 400}%
\special{fp}%
\special{pa 1600 400}%
\special{pa 1600 600}%
\special{fp}%
\special{pa 1600 600}%
\special{pa 400 600}%
\special{fp}%
\special{pa 400 400}%
\special{pa 400 600}%
\special{fp}%
\special{pa 600 600}%
\special{pa 600 400}%
\special{fp}%
\special{pa 800 400}%
\special{pa 800 600}%
\special{fp}%
\special{pa 1000 600}%
\special{pa 1000 400}%
\special{fp}%
\special{pa 1200 400}%
\special{pa 1200 600}%
\special{fp}%
\special{pa 1400 600}%
\special{pa 1400 400}%
\special{fp}%
%
\special{pn 8}%
\special{pa 1800 400}%
\special{pa 2800 400}%
\special{fp}%
\special{pa 2800 400}%
\special{pa 2800 600}%
\special{fp}%
\special{pa 2800 600}%
\special{pa 1800 600}%
\special{fp}%
\special{pa 1800 400}%
\special{pa 1800 800}%
\special{fp}%
\special{pa 1800 800}%
\special{pa 2000 800}%
\special{fp}%
\special{pa 2000 800}%
\special{pa 2000 400}%
\special{fp}%
\special{pa 2200 400}%
\special{pa 2200 600}%
\special{fp}%
\special{pa 2400 600}%
\special{pa 2400 400}%
\special{fp}%
\special{pa 2600 400}%
\special{pa 2600 600}%
\special{fp}%
%
\special{pn 8}%
\special{pa 3000 400}%
\special{pa 3000 800}%
\special{fp}%
\special{pa 3000 800}%
\special{pa 3400 800}%
\special{fp}%
\special{pa 3400 800}%
\special{pa 3400 400}%
\special{fp}%
\special{pa 3000 400}%
\special{pa 3800 400}%
\special{fp}%
\special{pa 3800 400}%
\special{pa 3800 600}%
\special{fp}%
\special{pa 3800 600}%
\special{pa 3000 600}%
\special{fp}%
\special{pa 3200 400}%
\special{pa 3200 800}%
\special{fp}%
\special{pa 3600 600}%
\special{pa 3600 400}%
\special{fp}%
%
\special{pn 8}%
\special{pa 4000 400}%
\special{pa 4000 1000}%
\special{fp}%
\special{pa 4000 1000}%
\special{pa 4200 1000}%
\special{fp}%
\special{pa 4200 1000}%
\special{pa 4200 400}%
\special{fp}%
\special{pa 4000 400}%
\special{pa 4600 400}%
\special{fp}%
\special{pa 4600 400}%
\special{pa 4600 600}%
\special{fp}%
\special{pa 4600 600}%
\special{pa 4000 600}%
\special{fp}%
\special{pa 4000 800}%
\special{pa 4400 800}%
\special{fp}%
\special{pa 4400 800}%
\special{pa 4400 400}%
\special{fp}%

\put(4.4000,-5.5000){\makebox(0,0)[lb]{$6$}}%

\put(6.4000,-5.5000){\makebox(0,0)[lb]{$5$}}%

\put(8.4000,-5.5000){\makebox(0,0)[lb]{$4$}}%

\put(10.4000,-5.5000){\makebox(0,0)[lb]{$3$}}%

\put(12.4000,-5.5000){\makebox(0,0)[lb]{$2$}}%

\put(14.4000,-5.5000){\makebox(0,0)[lb]{$1$}}%

\put(18.4000,-5.5000){\makebox(0,0)[lb]{$6$}}%

\put(30.4000,-5.5000){\makebox(0,0)[lb]{$5$}}%

\put(32.4000,-5.5000){\makebox(0,0)[lb]{$4$}}%

\put(40.4000,-5.5000){\makebox(0,0)[lb]{$5$}}%
\end{picture}%
\\
The number of odd numbers equals the number of even numbers. 
\begin{proof}
In Theorem 3.2, the strict partitions they have same arm length and different parity length are pair. 
Recall that $h_{1j}(\lambda ) = a_{1j}(\lambda ) + \ell (\lambda )-1$. 
Therefore the parity of their hook length are different. 
And the leftovers are divisors of $n$. 
When $n$ equals $2(2m+1)$, the number of odd divisors of $n$ equals the number of even divisors of $n$. 
\end{proof}

\section{Generalizations}

\begin{thm} For any positive integers $k,m,n$, 
\begin{eqnarray*}
&\sharp \{ (\lambda ,  i_1 ,  \ldots , i_{m})\ |\ \lambda \in {\mathcal {SP} }(n), 1\leq i_1 < \ldots < i_m \leq \lambda_{\ell (\lambda )} , a_{1,i_m}(\lambda )=k , \ell (\lambda )\textrm {is odd} \} \\
-&\sharp \{ (\lambda , i_1 , \ldots , i_{m} )\ |\ \lambda \in {\mathcal {SP} }(n), 1\leq i_1 < \ldots < i_m \leq \lambda_{\ell (\lambda )}, a_{1,i_m}(\lambda )=k , \ell (\lambda )\textrm {is even} \} \\
=&
\sharp \{ (\lambda , t_1, \ldots , t_{m- c(\lambda )} ) \ |\ \lambda \in {\mathcal {P}}(n), 1\leq t_1 < \ldots <t_{m-c(\lambda )}<\ell (\lambda ), \lambda _1 = k, \lambda_{t_i}= \lambda_{t_i +1} \} . 
\end{eqnarray*}

\end{thm}
\noindent
When $m=1$, this identity is equivalent to Theorem 3.2. \\

\bigskip
\noindent
\textbf{Example. }@For $n=5, m=2$, 
we draw the same figure as Figure 2. \\
\\
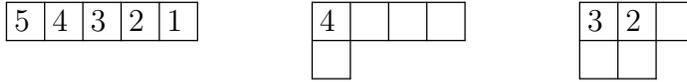
\textbf{Figure 4.}\\
\unitlength 0.1in
\begin{picture}(38.00,10.00)(4.00,-14.00)
%
\special{pn 8}%
\special{pa 400 400}%
\special{pa 400 400}%
\special{fp}%
%
\special{pn 8}%
\special{pa 600 600}%
\special{pa 1600 600}%
\special{fp}%
\special{pa 1600 600}%
\special{pa 1600 800}%
\special{fp}%
\special{pa 1600 800}%
\special{pa 600 800}%
\special{fp}%
\special{pa 600 800}%
\special{pa 600 600}%
\special{fp}%
%
\special{pn 8}%
\special{pa 2200 600}%
\special{pa 3000 600}%
\special{fp}%
\special{pa 3000 600}%
\special{pa 3000 800}%
\special{fp}%
\special{pa 3000 800}%
\special{pa 2200 800}%
\special{fp}%
\special{pa 2200 600}%
\special{pa 2200 1000}%
\special{fp}%
\special{pa 2200 1000}%
\special{pa 2400 1000}%
\special{fp}%
\special{pa 2400 1000}%
\special{pa 2400 600}%
\special{fp}%
\special{pa 2600 600}%
\special{pa 2600 800}%
\special{fp}%
\special{pa 2800 800}%
\special{pa 2800 600}%
\special{fp}%
\special{pa 1400 600}%
\special{pa 1400 800}%
\special{fp}%
\special{pa 1200 800}%
\special{pa 1200 600}%
\special{fp}%
\special{pa 1000 600}%
\special{pa 1000 800}%
\special{fp}%
\special{pa 800 800}%
\special{pa 800 600}%
\special{fp}%
%
\special{pn 8}%
\special{pa 3600 600}%
\special{pa 4200 600}%
\special{fp}%
\special{pa 4200 600}%
\special{pa 4200 800}%
\special{fp}%
\special{pa 4200 800}%
\special{pa 3600 800}%
\special{fp}%
\special{pa 3600 600}%
\special{pa 3600 1000}%
\special{fp}%
\special{pa 3600 1000}%
\special{pa 4000 1000}%
\special{fp}%
\special{pa 4000 1000}%
\special{pa 4000 600}%
\special{fp}%
\special{pa 3800 600}%
\special{pa 3800 1000}%
\special{fp}%

\put(6.4000,-7.5000){\makebox(0,0)[lb]{$5$}}%

\put(8.4000,-7.5000){\makebox(0,0)[lb]{$4$}}%

\put(10.4000,-7.5000){\makebox(0,0)[lb]{$3$}}%

\put(12.4000,-7.5000){\makebox(0,0)[lb]{$2$}}%

\put(14.4000,-7.5000){\makebox(0,0)[lb]{$1$}}%

\put(22.4000,-7.5000){\makebox(0,0)[lb]{$4$}}%

\put(36.4000,-7.5000){\makebox(0,0)[lb]{$3$}}%

\put(38.4000,-7.5000){\makebox(0,0)[lb]{$2$}}%

\end{picture}%
\\
We count the pairs of arm lengths in the same partition $\lambda $ with sign $(-1)^{\ell (\lambda )}$. 
The pairs with positive sign are $(5,1)$, $(4,1)$, $(3,1)$, $(2,1)$, $(5,2)$, $(4,2)$, $(3,2)$, $(5,3)$, $(4,3)$, $(5,4)$. 
The pair with negative sign is $(3,2)$. 
On the other hand, the Young diagrams of size $5$ they are made by concatenating $2$ rectangles are, \\
\\
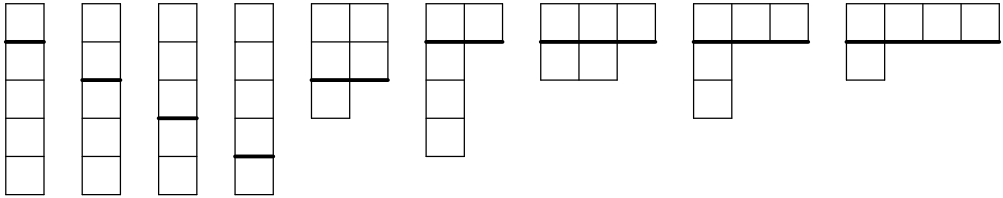
\textbf{Figure 5.}\\

\bigskip

\unitlength 0.1in
\begin{picture}(52.00,10.00)(4.00,-14.00)

\special{pn 8}%
\special{pa 400 400}%
\special{pa 400 1400}%
\special{fp}%
\special{pa 400 1400}%
\special{pa 600 1400}%
\special{fp}%
\special{pa 600 1400}%
\special{pa 600 400}%
\special{fp}%
\special{pa 600 400}%
\special{pa 400 400}%
\special{fp}%
\special{pa 400 800}%
\special{pa 600 800}%
\special{fp}%
\special{pa 600 1000}%
\special{pa 400 1000}%
\special{fp}%
\special{pa 400 1200}%
\special{pa 600 1200}%
\special{fp}%
\special{pa 800 400}%
\special{pa 800 1400}%
\special{fp}%
\special{pa 800 1400}%
\special{pa 1000 1400}%
\special{fp}%
\special{pa 1000 1400}%
\special{pa 1000 400}%
\special{fp}%
\special{pa 1000 400}%
\special{pa 800 400}%
\special{fp}%
\special{pa 800 600}%
\special{pa 1000 600}%
\special{fp}%
\special{pa 1000 1000}%
\special{pa 800 1000}%
\special{fp}%
\special{pa 800 1200}%
\special{pa 1000 1200}%
\special{fp}%
\special{pa 1200 1400}%
\special{pa 1200 400}%
\special{fp}%
\special{pa 1200 400}%
\special{pa 1400 400}%
\special{fp}%
\special{pa 1400 400}%
\special{pa 1400 1400}%
\special{fp}%
\special{pa 1400 1400}%
\special{pa 1200 1400}%
\special{fp}%
\special{pa 1200 600}%
\special{pa 1400 600}%
\special{fp}%
\special{pa 1400 800}%
\special{pa 1200 800}%
\special{fp}%
\special{pa 1200 1200}%
\special{pa 1400 1200}%
\special{fp}%
\special{pa 1600 600}%
\special{pa 1600 600}%
\special{fp}%

\special{pn 8}%
\special{pa 1600 400}%
\special{pa 1600 1400}%
\special{fp}%
\special{pa 1600 1400}%
\special{pa 1800 1400}%
\special{fp}%
\special{pa 1800 1400}%
\special{pa 1800 400}%
\special{fp}%
\special{pa 1800 400}%
\special{pa 1600 400}%
\special{fp}%
\special{pa 1600 600}%
\special{pa 1800 600}%
\special{fp}%
\special{pa 1800 800}%
\special{pa 1600 800}%
\special{fp}%
\special{pa 1600 1000}%
\special{pa 1800 1000}%
\special{fp}%
\special{pa 2000 400}%
\special{pa 2400 400}%
\special{fp}%
\special{pa 2400 400}%
\special{pa 2400 800}%
\special{fp}%
\special{pa 2000 400}%
\special{pa 2000 1000}%
\special{fp}%
\special{pa 2000 1000}%
\special{pa 2200 1000}%
\special{fp}%
\special{pa 2200 1000}%
\special{pa 2200 400}%
\special{fp}%
\special{pa 2000 600}%
\special{pa 2400 600}%
\special{fp}%
\special{pa 2600 400}%
\special{pa 3000 400}%
\special{fp}%
\special{pa 3000 400}%
\special{pa 3000 600}%
\special{fp}%
\special{pa 2600 400}%
\special{pa 2600 1200}%
\special{fp}%
\special{pa 2600 1200}%
\special{pa 2800 1200}%
\special{fp}%
\special{pa 2800 1200}%
\special{pa 2800 400}%
\special{fp}%
\special{pa 2600 800}%
\special{pa 2800 800}%
\special{fp}%
\special{pa 2800 1000}%
\special{pa 2600 1000}%
\special{fp}%
\special{pa 3200 400}%
\special{pa 3800 400}%
\special{fp}%
\special{pa 3200 400}%
\special{pa 3200 800}%
\special{fp}%
\special{pa 3200 800}%
\special{pa 3600 800}%
\special{fp}%
\special{pa 3600 800}%
\special{pa 3600 400}%
\special{fp}%
\special{pa 3400 400}%
\special{pa 3400 800}%
\special{fp}%
\special{pa 3800 400}%
\special{pa 3800 600}%
\special{fp}%

\special{pn 8}%
\special{pa 4000 400}%
\special{pa 4600 400}%
\special{fp}%
\special{pa 4600 400}%
\special{pa 4600 600}%
\special{fp}%
\special{pa 4000 400}%
\special{pa 4000 1000}%
\special{fp}%
\special{pa 4000 1000}%
\special{pa 4200 1000}%
\special{fp}%
\special{pa 4200 1000}%
\special{pa 4200 400}%
\special{fp}%
\special{pa 4400 400}%
\special{pa 4400 600}%
\special{fp}%
\special{pa 4200 800}%
\special{pa 4000 800}%
\special{fp}%
\special{pa 4800 400}%
\special{pa 5600 400}%
\special{fp}%
\special{pa 4800 400}%
\special{pa 4800 800}%
\special{fp}%
\special{pa 4800 800}%
\special{pa 5000 800}%
\special{fp}%
\special{pa 5000 800}%
\special{pa 5000 600}%
\special{fp}%
\special{pa 5600 400}%
\special{pa 5600 600}%
\special{fp}%
\special{pa 5400 600}%
\special{pa 5400 400}%
\special{fp}%
\special{pa 5200 400}%
\special{pa 5200 600}%
\special{fp}%
\special{pa 5000 600}%
\special{pa 5000 400}%
\special{fp}%

\special{pn 20}%
\special{pa 400 600}%
\special{pa 600 600}%
\special{fp}%
\special{pa 800 800}%
\special{pa 1000 800}%
\special{fp}%
\special{pa 1200 1000}%
\special{pa 1400 1000}%
\special{fp}%
\special{pa 1600 1200}%
\special{pa 1800 1200}%
\special{fp}%
\special{pa 2000 800}%
\special{pa 2400 800}%
\special{fp}%
\special{pa 2600 600}%
\special{pa 3000 600}%
\special{fp}%
\special{pa 3200 600}%
\special{pa 3800 600}%
\special{fp}%
\special{pa 4000 600}%
\special{pa 4600 600}%
\special{fp}%
\special{pa 4800 600}%
\special{pa 5600 600}%
\special{fp}%
\end{picture}%
\\

\bigskip

\noindent
There are $9$ such Young diagrams. 
Here we count the ways of concatenating rectangles. 
Hence, for example, the first $4$ diagrams must be thought of as the different ones. 
Therefore the number of pairs with positive sign minus the number of pairs with negative sign equals the number of Young diagrams made by concatenating $2$ rectangles. 
And more fix smaller number b, the number of pair $(a,b)$ with positive sign minus the number of pair $(a,b)$ with negative sign also equals 
the number of Young diagrams $\lambda$ made by concatenating $2$ rectangles that $\lambda _1$ equals $b$. 
For example, the number of the pairs that smaller number is $1$ and the number of Young diagrams that size of first line is $1$ are both $4$. 

\bigskip

\begin{proof}
When $\lambda \in {\mathcal {SP}}(n) $, $i_1 <\ldots i_m \leq \lambda _\ell $ are given, we make strict partitions $\lambda ^{(1)}, \ldots , \lambda ^{(m)}$ in the following procedure. 
First, we put $\lambda ^{(1)} $ with $\lambda $. 
When $\lambda ^{(h)}$ is determined, we put $\mu=\lambda ^{(h)} $ and $j_{h} =a_{1,i_{m+1-h } }(\lambda ^{(h ) } ) $. 
And we make new strict partition [$\mu $] by replacing the maximum parts $\mu_1$ by $\mu_1 - j_h$. 
This operation keeps strictness, since $\mu $ is strict in ${\textrm {mod}} k$. 
We repeat this operation until we get $a_{1, i_{m-h } }(\mu )\leq j_h $. 
We put $\lambda ^{(h+1)}$ with $\mu $ obtained in this way. 
And let $t_h$ be the number of times of the operations. 
The lengths of $\lambda ^{(i )}$-s are same as the length of $\lambda $. 
We consider that $j_m $ is $k$ of theorem 3.2, and pair up $\lambda ^{(m)}$-s by $\alpha _{j_m}$. 
Then there is only $1$ difference between the lengths of partitions of each pair. 
Leftovers are $\lambda ^{(m)}=(\lambda^{(m)}_1 )$ that $\lambda^{(m)}_1 $ is multiple of $j_m$. 
Since $j_1 \geq \ldots \geq j_m$, it corresponds with Young diagram that is made by concatenating rectangles $j_h\times {t_h} $, 
where $\displaystyle t_m=\frac{\lambda^{(m)}_1 }{j_m}$. 

\end{proof}

Taking the total sum over $k$ for Theorem 4.1, we have an identity which is shown in the following generating function. 
\begin{thm}$($Dilcher's identity 1$)$
For any positive integer $m$,  
\[
\displaystyle \sum_{k= 1}^{\infty }{(-1)^k \frac {q^{\frac {k(k+1) }{2 }
+(m-1)k }
 }{(q;q )_k (1-q^k)^m } }
=\sum_{j_1= 1 }^{\infty }{\frac {q^{j_1} }{1-q^{j_1} } 
\sum_{j_2 =1 }^{j_1 }{\frac {q^{j_2} }{1-q^{j_2} } 
\cdots
\sum_{j_m =1 }^{j_{m-1 } }{\frac {q^{j_m} }{1-q^{j_m} } 
}
}
}. 
\]
\end{thm}
\noindent
In a language of Young diagrams, this identity is equivalent to the following. 
\[
\displaystyle
\sum_{\lambda \in {\mathcal {SP} } }{(-1)^{\ell (\lambda)-1 }\binom {\lambda _{\ell (\lambda) }}{m} q^{|\lambda | } } =
\sum_{
\lambda \in {\mathcal {P} }
}
{\binom {\ell (\lambda )-c(\lambda )}{m-c(\lambda ) }q^{|\lambda |}
}. 
\]
And taking the sum over $k$ from $1$ to $k$ for Theorem 4.1, we have an identity which is shown in the following generating function. 

\begin{thm}$($Dilcher's identity 2$)$
For any positive integers $m$ and $t $, 
\[
\displaystyle \sum_{k= 1}^{t }{(-1)^k \frac {q^{\frac {k(k+1) }{2 }
+(m-1)k }
 }{(1-q^k)^m } 
\left [ 
\begin{array}{l} 
t \\
k 
\end{array}
\right ] 
}
=\sum_{j_1= 1 }^{t }{\frac {q^{j_1} }{1-q^{j_1} } 
\sum_{j_2 =1 }^{j_1 }{\frac {q^{j_2} }{1-q^{j_2} } 
\cdots
\sum_{j_m =1 }^{j_{m-1 } }{\frac {q^{j_m} }{1-q^{j_m} } 
}
}
}. 
\]
\end{thm}

\bigskip
\noindent
Theorem 4.1 is a generalization of Dilcher's identities. 
Which is the coefficient of $b^kq^n$ in the following.
\begin{thm}
For any positive integer $m$, 
\[
\displaystyle \sum_{k= 1}^{\infty }{(-1)^k \frac {b^kq^{\frac {k(k+1) }{2 }
+(m-1)k }
 }{(bq;q )_k (1-q^k)^m } }
=\sum_{j_1= 1 }^{\infty }{\frac {b^{j_1 }q^{j_1} }{1-q^{j_1} } 
\sum_{j_2 =1 }^{j_1 }{\frac {q^{j_2} }{1-q^{j_2} } 
\cdots
\sum_{j_m =1 }^{j_{m-1 } }{\frac {q^{j_m} }{1-q^{j_m} } 
}
}
}. 
\]
\end{thm}
\begin{thm}
For any positive integers $m$ and $t$, 
\[
\displaystyle \sum_{k= 1}^{t }{(-1)^k \frac {b^kq^{\frac {k(k+1) }{2 }
+(m-1)k }
 }{(1-q^k)^m } 
\left [ 
\begin{array}{l} 
t \\
k 
\end{array}
\right ] _{q,b}
}
=\sum_{j_1= 1 }^{t }{\frac {b^{j_1 }q^{j_1} }{1-q^{j_1} } 
\sum_{j_2 =1 }^{j_1 }{\frac {q^{j_2} }{1-q^{j_2} } 
\cdots
\sum_{j_m =1 }^{j_{m-1 } }{\frac {q^{j_m} }{1-q^{j_m} } 
}
}
},  
\]
where $\displaystyle \left [ 
\begin{array}{l} 
t \\
k 
\end{array}
\right ] _{q,b}$
is defined as
\[
\left [ 
\begin{array}{l} 
t \\
k 
\end{array}
\right ] _{q,b}
= 
\sum _{
\begin{array}{c}
\lambda \in {\mathcal {P }}\\
\lambda _1 \leq t-k , \ell (\lambda )\leq k
\end{array}
}
{b^{\lambda _1}q^{|\lambda |}}. 
\] 
\end{thm}
\noindent
When $b=1$, they are Dilcher's identities. 

\end{document}